\documentclass[11pt]{article}
\usepackage{amsmath, mathrsfs}
\usepackage{amssymb,amsmath,amsfonts}
\usepackage{hyperref}
\usepackage{color}
\newtheorem{theorem}{Theorem}[section]

\makeatletter \@addtoreset{equation}{section} \makeatother
  \textheight
21.5cm \textwidth 15.2cm \evensidemargin 0.2cm \oddsidemargin 0.2cm

\begin{document}

\bigskip

\begin{center}
{\Large {\bf A note on generalized \boldmath $S$-space-forms}}

\bigskip

Mukut Mani Tripathi

\bigskip
\end{center}

\noindent {\bf Abstract.} It is shown that, in a generalized $S$-space-form,
$F_{7}$ always equals $F_{8}$. \medskip

 \medskip


\section{Introduction\label{sect-intro}}

A $(2m+s)$-dimensional Riemannian manifold $(M,g)$ endowed with an $f$%
-structure $f$ (that is, a tensor field of type $(1,1)$ and rank $2m$
satisfying $f^{3}+f=0$ \cite{Y}) is said to be a {\em metric} $f${\em %
-manifold} if there exist $s$ global vector fields $\xi _{1},\dots ,\xi _{s}$
on $M$ (called {\em structure vector fields}) such that, if $\eta _{1},\dots
,\eta _{s}$ are the dual 1-forms of $\xi _{1},\dots ,\xi _{s}$, then
\[
f\xi _{\alpha }=0;\quad \eta _{\alpha }\circ f=0;\quad
f^{2}=-\,I+\sum_{\alpha =1}^{s}\eta _{\alpha }\otimes \xi _{\alpha }\,,
\]%
\[
g(X,Y)=g(fX,fY)+\sum_{\alpha =1}^{s}\eta _{\alpha }(X)\eta _{\alpha }(Y)
\]%
for all $X,Y\in {\cal X}(M)$ and $\alpha =1,\ldots ,s$. Let $F$ be the
fundamental $2$-form on $M$ defined by $F(X,Y)=g(X,fY)$ for any $X,Y\in
{\cal X}(M)$. A metric $f$-manifold is {\em normal} if the Nijenhuis tensor $%
\left[ f,f\right] $ of $f$ equals $-\,2\sum_{\alpha =1}^{s}\xi _{\alpha
}\otimes d\eta _{\alpha }$. A normal metric $f$-manifold is said to be an $S$%
{\em -manifold} if $F={\rm d}\eta _{\alpha }$ for each $\alpha \in \left\{
1,\dots ,s\right\} $, while a normal metric $f$-manifold is a $K${\em %
-manifold} if ${\rm d}F=0$. A $K$-manifold is a $C$-manifold if ${\rm d}\eta
_{\alpha }=0$ for each $\alpha \in \left\{1,\dots ,s\right\} $. For $s=0$, a
$K$-manifold is a Kaehler manifold and, for $s=1$, a $K$-manifold is a
quasi-Sasakian manifold, an $S$-manifold is a Sasakian manifold and a $C$%
-manifold is a cosymplectic manifold. For more details we refer to \cite%
{B,CFF-09-AG}. \medskip

A plane section $\Pi$ on a metric $f$-manifold $M$ is said to be an $f${\em %
-section} if it is determined by a unit vector $X$ orthogonal to all
structure vector fields and the unit vector field $fX$. The sectional
curvature of $\Pi $ is called an $f${\em -sectional curvature}. An $S$%
-manifold (resp., a $C$-manifold) is said to be an $S${\em -space-form}
(resp., a $C${\em -space-form}) if it has constant $f$-sectional curvature $c
$ and then, it is denoted by $M(c)$ \cite{KT}. \medskip

As a generalization of $S$-space-forms and $C$-space-forms, in
\cite{CFF-09-AG} the authors introduced the concept of a generalized
$S$-space-form equipped with two structure vector fields and studied
its basic geometric properties. In fact, a metric $f$-manifold
$\left( M,f,\xi _{1},\xi _{2},\eta _{1},\eta _{2},g\right) $ with
two structure vector fields $\xi _{1}$ and $\xi _{2}$ is called a
{\em generalized} $S${\em -space-form} \cite{CFF-09-AG} if there
exist differentiable functions $F_{1},\dots ,F_{8}$ on $M$ such that
the curvature tensor field $R$ of $M$ satisfies
\begin{eqnarray}
R(X,Y)Z &=&F_{1}\left\{ g(Y,Z)X-g(X,Z)Y\right\}  \nonumber \\
&&+ F_{2}\left\{ g(X,fZ)fY-g(Y,fZ)fX+2g(X,fY)fZ\right\}  \nonumber \\
&&+ F_{3}\left\{ \eta _{1}(X)\eta _{1}(Z)Y-\eta _{1}(Y)\eta
_{1}(Z)X+g(X,Z)\eta _{1}(Y)\xi _{1}-g(Y,Z)\eta _{1}(X)\xi
_{1}\right\}
\nonumber \\
&&+ F_{4}\left\{ \eta _{2}(X)\eta _{2}(Z)Y-\eta _{2}(Y)\eta
_{2}(Z)X+g(X,Z)\eta _{2}(Y)\xi _{2}-g(Y,Z)\eta _{2}(X)\xi
_{2}\right\}
\nonumber \\
&&+ F_{5}\left\{ \eta _{1}(X)\eta _{2}(Z)Y-\eta _{1}(Y)\eta
_{2}(Z)X+g(X,Z)\eta _{1}(Y)\xi _{2}-g(Y,Z)\eta _{1}(X)\xi
_{2}\right\}
\nonumber \\
&&+ F_{6}\left\{ \eta _{2}(X)\eta _{1}(Z)Y-\eta _{2}(Y)\eta
_{1}(Z)X+g(X,Z)\eta _{2}(Y)\xi _{1}-g(Y,Z)\eta _{2}(X)\xi
_{1}\right\}
\nonumber \\
&&+ F_{7}\left\{ \eta _{1}(X)\eta _{2}(Y)\eta _{2}(Z)\xi _{1}-\eta
_{2}(X)\eta _{1}(Y)\eta _{2}(Z)\xi _{1}\right\}  \nonumber \\
&&+ F_{8}\left\{ \eta _{2}(X)\eta _{1}(Y)\eta _{1}(Z)\xi _{2}-\eta
_{1}(X)\eta _{2}(Y)\eta _{1}(Z)\xi _{2}\right\}   \label{eq-gssf}
\end{eqnarray}%
for all $X,Y,Z\in {\cal X}(M)$. \medskip

Here, we prove the following

\begin{theorem}
\label{th-GSSF-F7=F8} Let $( M,f,\xi _{1},\xi _{2},\eta _{1},\eta
_{2},g) $ be a generalized\ $S$-space-form. Then $F_{7}=F_{8}$.
\end{theorem}

\section{Proof of Theorem~\protect\ref{th-GSSF-F7=F8}\label%
{sect-proof-GSSF-F7=F8}}

From (\ref{eq-gssf}) we obtain
\begin{equation}
R\left( \xi _{1},\xi _{2}\right) \xi _{1}=-\,F_{1}\xi _{2}+F_{3}\xi
_{2}+F_{4}\xi _{2}-F_{8}\xi _{2}\,,  \label{eq-1}
\end{equation}%
\begin{equation}
R\left( \xi _{1},\xi _{2}\right) \xi _{2}=F_{1}\xi _{1}-F_{3}\xi
_{1}-F_{4}\xi _{1}+F_{7}\xi _{1}.  \label{eq-2}
\end{equation}%
Using (\ref{eq-1}) and (\ref{eq-2}) in
\[
g\left( R\left( \xi _{1},\xi _{2}\right) \xi _{1}\,,\,\xi
_{2}\right) +g\left( R\left( \xi _{1},\xi _{2}\right) \xi
_{2}\,,\,\xi _{1}\right)=0
\]%
we get $F_{7}=F_{8}$. $\blacksquare $

\noindent Department of Mathematics

\noindent Faculty of Science

\noindent Banaras Hindu University

\noindent Varanasi 221005, India

\noindent E-mail: {\tt mmtripathi66@yahoo.com}

\newpage

\end{document}